\documentclass[runningheads]{llncs}
\usepackage{algorithmic}
\usepackage{algorithm}
\usepackage{mathrsfs}
\usepackage{latexsym}
\usepackage{amsmath}
\usepackage{amssymb}
\usepackage{color}
\usepackage{graphicx}
\usepackage{enumerate}
\usepackage{changebar}
\usepackage{url}
\usepackage{cases}
\usepackage{multirow}

\usepackage{rotating}

\makeatletter
\let\markeq\ltx@label
\makeatother

\usepackage{url}

\newcommand{\keywords}[1]{\par\addvspace\baselineskip
\noindent\keywordname\enspace\ignorespaces#1}

\begin{document}

\mainmatter  % start of an individual contribution

\title{Using a conic bundle method to accelerate both phases of a quadratic convex reformulation}

\titlerunning{Using Conic Bundle to accelerate a quadratic convex reformulation}

\author{Alain Billionnet\inst{1}, Sourour Elloumi\inst{1}, Am\'elie Lambert\inst{2}, Angelika Wiegele\inst{3}}
\institute{ENSIIE-CEDRIC, 1 square de la r\'esistance FR-91025 Evry\\
  \email{alain.billionnet@ensiie.fr, sourour.elloumi@ensiie.fr} 
  \and CNAM-CEDRIC, 292 Rue St Martin FR-75141 Paris Cedex 03\\
  \email{amelie.lambert@cnam.fr} 
  \and Institut f\"ur Mathematik, Alpen-Adria-Universit\"at Klagenfurt\\
  \email{angelika.wiegele@aau.at}}

\authorrunning{Alain Billionnet, Sourour Elloumi, Am\'elie Lambert, Angelika Wiegele}

\maketitle
 
\begin{abstract}
  We present algorithm \texttt{MIQCR-CB} that is an advancement of method
  \texttt{MIQCR}~(Billionnet, Elloumi and Lambert, 2012). \texttt{MIQCR} is a method for solving mixed-integer quadratic
  programs and works in two phases: the first phase determines an
  equivalent quadratic formulation with a convex objective function by solving a
  semidefinite problem $(SDP)$, and, in the second phase, the equivalent formulation
  is solved by a standard solver. As the reformulation relies on the solution
  of a large-scale semidefinite program, it is not tractable by existing semidefinite solvers, already for medium
  sized problems. To surmount this
  difficulty, we present in \texttt{MIQCR-CB} a subgradient
  algorithm within a Lagrangian duality framework for solving $(SDP)$ that substantially speeds up the first phase. Moreover, this algorithm leads to a reformulated problem of smaller size than the one obtained by the original \texttt{MIQCR} method which results in a shorter time for solving the second phase.
 We present extensive computational results to show the efficiency of our algorithm. First, we apply \texttt{MIQCR-CB} to the $k$-cluster problem that can be
 formulated by a binary quadratic program. As an illustration of the
 efficiency of our new algorithm, for instances of size 80 and of density 25$\%$,
 \texttt{MIQCR-CB} is on average 78~times faster for Phase~1 and 24 times
 faster for Phase~2 than the original \texttt{MIQCR}. We also compare \texttt{MIQCR-CB} with \texttt{QCR}~(Billionnet, Elloumi and Plateau, 2009) and with \texttt{BiqCrunch}~(Krislock, Malick and
 Roupin, 2013) two methods devoted to binary quadratic programming.  We show that
 \texttt{MIQCR-CB} is able to solve most of the $225$ considered instances within
 $3$ hours of cpu time. We also present experiments on two classes of general integer instances where we compare \texttt{MIQCR-CB} with \texttt{MIQCR}, \texttt{Couenne} and \texttt{Cplex12.6}. We demonstrate the significant improvement over the original \texttt{MIQCR} approach. Finally, we show that \texttt{MIQCR-CB} is able to solve almost all of the considered instances while \texttt{Couenne} and \texttt{Cplex12.6} are not able to solve half out of them.

\end{abstract}

\section{Introduction}

We present an algorithm that accelerates the computation time of method \texttt{MIQCR (Mixed-Integer Quadratic Convex Reformulation)}~\cite{BEL12}. This method is an exact solution algorithm for quadratic programs having mixed-integer variables. This obviously includes the class of quadratic problems having pure integer variables such as $(QP)$:
\begin{numcases}{(QP)}
\max  f(x)= \displaystyle{\sum_{i=1}^n} \displaystyle{\sum_{j=i}^n}  q_{ij} x_i x_j + \sum_{i=1}^n c_ix_i \nonumber \\
\mbox{s.t.} \nonumber \\
\quad  \displaystyle{\sum_{i=1}^n} a_{ri}x_i = b_r  & $1 \le r \le m$ \label{equality}\\
\quad  0 \leq x_i \leq u_i &  $i \in I$ \label{bound}\\
\quad x_i \in \mathbb{N}_{0}   &  $i \in I$ \label{integer} 
\end{numcases}

\noindent where $I=\{1,\ldots,n\}$.

$(QP)$ belongs to the class of Mixed-Integer Non-Linear Programs (MINLP). It includes the case of linear inequalities, since any problem with inequalities can be rewritten as $(QP)$ by introducing non-negative slack variables. Conventional approaches to solve MINLP are based on global optimization techniques~\cite{AudHanJauSav00,BLL09,Floudas00,LibMac06,SahTaw05,TawSah02}. Software is available to solve this large class of problems that includes $(QP)$, see for instance~\cite{couenne,baron} . The solver \texttt{Couenne}~\cite{couenne} uses linear relaxations within a spatial branch-and-bound algorithm~\cite{Smi96,SmiPan97,SmiPan99}, together with heuristics for finding feasible solutions. Using a branch-and-bound framework based on convex relaxations~\cite{BliBon13}, the recent implementation of \texttt{Cplex12.6}~\cite{cplex126} also handles $(QP)$.

Many applications in operations research and industrial engineering involve general integer variables in their formulation. Some of these applications can be formulated as $(QP)$. Methods are available for solving particular cases of $(QP)$. If the objective function is linear, we refer to Mixed-Integer Linear Programming (MILP), which is still $\mathcal{NP}$-Hard, but for which a large variety of methods are well developed. If we assume the objective function to be convex, there also exist fairly efficient solvers~\cite{bonmin,cplex125}. 

\texttt{MIQCR} is an algorithm in two phases: the first phase computes an equivalent quadratic convex formulation of the initial problem by solving a large semidefinite problem, and in the second phase the reformulated problem is solved by a standard solver. Due to its size, the solution of the semidefinite problem of Phase~1 often constitutes the bottleneck of this method. However, once the equivalent formulation is computed, solving the obtained reformulated program is practical, since the continuous relaxation bound of the reformulation is tight. Hence, to handle larger instances method \texttt{MIQCR} needs an appropriate algorithm to solve Phase~1, while Phase~2 can still be handled by a standard solver. Thus, our first contribution in this paper lies in a different algorithm to solve Phase~1 of \texttt{MIQCR}. 
We first introduce a subgradient algorithm within a Lagrangian duality framework for solving $(SDP)$ approximately following the procedure introduced in~\cite{fischeretal}. Then, we parameterize our algorithm obtaining a dual heuristic for solving Phase~1 of the original \texttt{MIQCR} method. With this new algorithm, the time for computing Phase~1 significantly decreases and allows us to handle large-scale instances. Moreover, we obtain also a speed up of Phase~2 since by construction the equivalent formulation computed with the new algorithm is smaller than the reformulation obtained in the original \texttt{MIQCR} method. Hence, we can claim that our new algorithm \texttt{MIQCR-CB} is a general exact solution method for mixed-integer (or binary) quadratic programs of large size.

To illustrate our algorithm we apply it to the $k$-cluster problem. This problem can be formulated by an equality constrained binary quadratic program, the subclass of $(QP)$ obtained by setting all the upper bounds of the integer variables to one. The reason for choosing  a binary quadratic program is to demonstrate the impact of the new ideas compared to method~\texttt{QCR (Quadratic Convex Reformulation)}~\cite{BEP09} from both theoretical and experimental point of view. Recall that \texttt{MIQCR} exploits the ideas behind \texttt{QCR} and widens the applicability to the general mixed-integer case. As \texttt{MIQCR}, \texttt{QCR}~is also a method in two phases, where the computation of the equivalent formulation requires the solution of a semidefinite problem. The advantage of \texttt{QCR} over \texttt{MIQCR} lies in the short time required to compute the quadratic convex equivalent formulation (i.e. to solve the associated semidefinite problem). However, for large instances, this method is limited by the weakness of its bound. Hence, methods \texttt{QCR} and \texttt{MIQCR} do not have the same bottleneck: \texttt{QCR} is limited by the time needed for solving Phase~2 due to the weakness of its bound, and \texttt{MIQCR} is limited by the solution time of the huge semidefinite problem of Phase~1.
In the experiments the solution technique of Phase~1 of \texttt{MIQCR-CB} turns out to be almost as fast as Phase~1 in \texttt{QCR}, while the computed bound is as tight as in \texttt{MIQCR}. 

We compare experimentally \texttt{MIQCR-CB} with \texttt{MIQCR} and \texttt{QCR}. We also compare \texttt{MIQCR-CB} with the recent approach of Krislock, Malick and Roupin~\cite{KMR13} called \texttt{BiqCrunch}, on $225$ instances of the $k$-cluster problem with up to $160$ variables. Algorithm ~\texttt{BiqCrunch} is developed for solving binary quadratic programs. It consists of the branch-and-bound framework BOB~\cite{bob} using semidefinite programming bounds~\cite{Mal07}. We show that our approach is comparable with~\texttt{BiqCrunch} for solving instances with up to $120$ variables. \texttt{BiqCrunch} is slightly faster on larger instances, but, the limit of both algorithms (\texttt{BiqCrunch} and \texttt{MIQCR-CB}) lies in instances having 140 to 160 variables.

We also draw comparisons for general integer quadratic problems. First, we test our algorithm on Equality Integer Quadratic Problems (EIQP) and we show that \texttt{MIQCR-CB} is about 3 times faster than \texttt{MIQCR}. We also compare \texttt{MIQCR-CB} with the solvers \texttt{Couenne}~\cite{couenne} and \texttt{Cplex12.6}~\cite{cplex126} on an integer problem of equipartition that can be seen as the extension of a classical binary combinatorial optimization problem to the general integer case. For the considered instances,  \texttt{MIQCR-CB} is able to solve $27$ instances over the $40$ presented while the solvers \texttt{Cplex12.6} and \texttt{Couenne} solve only 3 and 8 instances respectively, within one hour of cpu time.

 The paper is organized as follows. In Section~\ref{sec:miqcr}, we describe algorithm \texttt{MIQCR} and discuss its limitations. Section~\ref{sec:solving} presents our improved algorithm \texttt{MIQCR-CB}. In Section~\ref{sec:comp}, we state the formal definition of the $k$-cluster problem and we present extensive computational results for solving  it, and in Section~\ref{sec:general} we show experimentally that our algorithm is efficient for the general integer case. Section~\ref{sec:conclusion} draws a conclusion.

\section{Recall of method \texttt{MIQCR}~\cite{BEL12}}
\label{sec:miqcr}

 In the following we describe method \texttt{MIQCR} (Mixed-Integer Quadratic Convex Reformulation)~\cite{BEL12}. In Phase~1 of \texttt{MIQCR} we aim to find an equivalent formulation of $(QP)$ having a concave objective function. Three parameters, $\alpha$, $\lambda$ and $\beta$, are obtained from the dual solution of a semidefinite relaxation of $(QP)$ that we call $(SDP)$. Using these parameters, we construct an equivalent problem $(QP_{\alpha, \lambda, \beta})$ having a concave objective function. The second phase consists of solving $(QP_{\alpha,\lambda,  \beta})$ by a standard solver. In detail this works as follows.\\

\noindent {\bf Phase~1: Constructing an equivalent formulation}\\
Consider the following semidefinite relaxation of $(QP)$.
\begin{numcases}{(SDP)}
\max f(X,x) = \displaystyle{\sum_{i=1}^n} \displaystyle{\sum_{j=i+1}^n}  q_ {ij} X_{ij} + \displaystyle{\sum_{i=1}^n}c_ix_i \nonumber  \\
\mbox{s.t.} \nonumber \\
  \displaystyle{\sum_{i=1}^n} a_{ri}x_i = b_r  & $1 \le r \le m$ (\ref{equality})  \nonumber\\ 
 \displaystyle{\sum_{r=1}^m} (\displaystyle{\sum_{i=1}^N} ( \displaystyle{\sum_{j=1}^N}a_{ri}a_{rj} X_{ij} - 2 a_{ri}b_rx_i)) = - \displaystyle{\sum_{r=1}^m}  b_r^2 & \label{c1}\\
 - X_{ii} +x_i  \leq 0 & $i \in I$  \label{i1}  \\ 
 - X_{ii} + 2u_ix_i - u_i^2  \leq 0 & $i \in I$  \label{i2}  \\
 X_{ii}  - u_ix_i  \leq 0& $i \in I$  \label{i3}  \\
 X_{ij} - u_jx_i \leq 0  & $(i,j) \in I^2, i< j $  \label{ij1}  \\ 
 X_{ij} -  u_ix_j \leq 0  & $(i,j) \in I^2, i< j $ \label{ij2}  \\ 
 - X_{ij} + u_jx_i + u_ix_j - u_iu_j\leq 0  &$(i,j) \in I^2, i< j $ \label{ij3}  \\ 
  - X_{ij} \leq 0  & $(i,j) \in I^2, i<j$  \label{ij4}  \\ 
  \left ( \begin{array}{ll}
   1 & x^T \\
   x & X 
     \end{array}\right ) \succeq 0  \label{sdp_1}  \\
 x \in \mathbb{R}^n,\;  X \in \mathcal{S}_n &   \label{sdp_2} 
\end{numcases}
where $\mathcal{S}_n$ is the space of symmetric matrices of order $n$. Constraint~(\ref{c1}) of $(SDP)$ is obtained by squaring and then summing up the $m$ equations~(\ref{equality}), and replacing each product $x_ix_j$ by $X_{ij}$. Constraints~(\ref{i1}) arise from $x_i^2 \geq x_i$ which is true if $x_i$ is a general integer variable, and  Constraints~(\ref{i2})-(\ref{ij4}) are the so-called McCormick inequalities that tighten the formulation~\cite{Cor76}. $(SDP)$ is known as the "SDP + RLT" semidefinite relaxation~\cite{Ans09}. \\

 Let $(\alpha,\lambda,\beta)$ be a dual optimal solution of $(SDP)$ where:
\begin{itemize}
\item  $\alpha \, \in\,  \mathbb{R}$ is the dual variable associated with Constraint~(\ref{c1}), 
\item $\lambda_i \, \in \, \mathbb{R} $ with $\lambda_i = - \lambda^1_i -  \lambda^2_i +   \lambda^3_i$ for any $i \in I$, where $\lambda^1_i, \lambda^2_i,\lambda^3_i$, are the non-negative dual variables associated with Constraints~(\ref{i1}),(\ref{i2}),(\ref{i3}), respectively,
\item $\beta_{ij} \, \in \, \mathbb{R}$ with $\beta_{ij}= \beta^1_{ij} + \beta^2_{ij} - \beta^3_{ij}- \beta^4_{ij}$ for any $(i,j)$ with $1\le i<j \le n$, where $\beta^1_{ij}\,\beta^2_{ij},\beta^3_{ij},\beta^4_{ij}$, are the non-negative dual variables associated with Constraints~(\ref{ij1}),~(\ref{ij2}),~(\ref{ij3}),~(\ref{ij4}), respectively. 
\end{itemize}
We introduce the following reformulated function:
%\begin{scriptsize}
\begin{eqnarray}
\lefteqn{f_{\alpha,\lambda,\beta} (x,y) =} \label{eq:f}\\
&& f(x) + \alpha
 \displaystyle{\sum_{r=1}^m} (\displaystyle{\sum_{i=1}^N}a_{ri}x_i - b)^2 + \displaystyle{\sum_{i=1: \lambda_i \neq 0}^n}  \lambda_i (y_{ii} - x^2_i) + \displaystyle{\sum_{i=1 }^n} \displaystyle{\sum_{j=i+1 : \beta_{ij} \neq 0}^n}  \beta_{ij} (y_{ij} - x_ix_j) \nonumber
\end{eqnarray}
%\end{scriptsize}
It is clear that $f_{\alpha,\lambda,\beta} (x,y) = f(x)$ if $(x,y)$ satisfies Constraints~(\ref{equality}), and $y_{ij} = x_ix_j$ for any $(i,j) \, \in \, I^2$. We use the linearization of the equality $y_{ij} = x_ix_j$ introduced in~\cite{BEL12} and define the set $S_{xyzt}$ containing the quadruplets $(x,y,z,t)$ satisfying the following conditions:

\begin{numcases}{}
 x_i  = \sum_{k=0}^{\lfloor \log (u_i)\rfloor} 2^k t_{ik} & $ i \in I $   \label{dec_x}   \\  
 y_{ij}  = \sum_{k=0}^{\lfloor \log(u_i)\rfloor} 2^k z_{ijk}& $ (i,j) \in I^2   $ \label{dec_y}    \\  
 z_{ijk}  \leq u_jt_{ik}  &$(i,k) \in E,\,\,j \in I $    \label{z1} \\
  z_{ijk}  \leq x_j  &$(i,k) \in E,\,\,j \in I   $   \label{z2}   \\
  z_{ijk}  \geq x_j - u_j(1 - t_{ik})  &$(i,k) \in E,\,\,j \in I $  \label{z3}\\
   z_{ijk}  \geq 0  &$(i,k) \in E,\,\,j \in I  $   \label{z4}   \\ 
  t_{ik}   \in \{0, 1\} &$ (i,k) \in E   $ \label{binary} \\
   y_{ii} \geq x_i&$ i \in I  ,\,\,\lambda_i \neq 0  $   \label{iy1}\\
   y_{ii} \geq 2u_ix_i - u_i^2&$ i \in I ,\,\,\lambda_i \neq 0  $   \label{iy2}\\
   y_{ii} \leq u_ix_i&$ i \in I ,\,\,\lambda_i \neq 0  $   \label{iy3}\\
   y_{ij} =y_{ji}& $ (i,j) \in I^2,\,\, i< j ,\,\,\beta_{ij} \neq 0   $  \label{1y2}   \\  
   y_{ij}  \leq u_j x_i& $(i,j) \in I^2 ,\,\, i< j ,\,\,\beta_{ij} \neq 0$    \label{1y3}\\
   y_{ij}  \leq u_ix_j &$(i,j) \in I^2,\,\, i< j ,\,\,\beta_{ij} \neq 0  $   \label{1y4}  \\  
  y_{ij}  \geq u_j x_i + u_i x_j-u_iu_j& $(i,j) \in I^2 ,\,\, i< j ,\,\,\beta_{ij} \neq 0$    \label{1y5}\\
   y_{ij}  \geq 0 &$(i,j) \in I^2,\,\, i<j ,\,\,\beta_{ij} \neq 0 $   \label{1y6}    
\end{numcases}
 \noindent with  $E=\{(i,k)\colon\, i\in I,\; k=0,\ldots \lfloor log(u_i) \rfloor \}$. \\

This allows us to reformulate $(QP)$ as problem $(QP_{\alpha, \lambda,  \beta})$, where $f_{\alpha,\lambda,\beta} (x,y)$ is a concave function and all constraints are linear:

\begin{numcases}{(QP_{\alpha, \lambda, \beta})}
\max  f_{\alpha,\lambda,\beta} (x,y) \nonumber \\
\mbox{s.t.} \nonumber \\
\quad  \displaystyle{\sum_{i=1}^n} a_{ri}x_i = b_r  & $1 \le r \le m$  (\ref{equality}) \nonumber\\
(x,y,z,t) \in S_{xyzt} \nonumber
\end{numcases}

\noindent {\bf Phase~2: Solving the reformulated problem $(QP_{\alpha,\lambda,\beta})$}\\
 $(QP_{\alpha,\lambda,  \beta})$ has a concave quadratic objective function and all its constraints are linear, thus computing the optimal value of the continuous relaxation of $(QP_{\alpha,\lambda, \beta})$ can be done in polynomial time. A general-purpose MIQP solver consisting of a branch-and-bound framework based on the continuous relaxation can solve $(QP_{\alpha,\lambda,  \beta})$.

Several properties can be deduced from the general mixed-integer case presented in~\cite{BEL12} which are the following.
\begin{itemize}
\item Problem $(QP_{\alpha,\lambda,  \beta})$ is an equivalent formulation of $(QP)$ and the two problems have therefore the same optimal values. 
\item The optimal solution value of the continuous relaxation of $(QP_{\alpha,\lambda,  \beta})$ is equal to the optimal solution value of $(SDP)$. Moreover, parameters $\alpha$, $\lambda$, and $\beta$ computed as described above, provide the tightest convex equivalent formulation of $(QP)$ in this reformulation scheme. 
\item From any feasible dual solution $(\bar{\alpha},\bar{\lambda},\bar{\beta})$ of $(SDP)$, the associated function $f_{\bar{\alpha},\bar{\lambda},\bar{\beta}}(x,y)$ is concave and thus our equivalent formulation is valid even if we are not able to solve $(SDP)$ to optimality.\\
\end{itemize}

\begin{property} From any $(\alpha,\lambda,  \beta)$ where $f_{\alpha,\lambda,\beta} (x,y)$ is not a concave function, we are always able to build a concave function $f_{\alpha,\lambda,\beta} (x,y)$ by taking $(\alpha,\lambda + \epsilon, \beta)$, where $\epsilon$ is the biggest eigenvalue of the Hessian matrix of function $f_{\alpha,\lambda,\beta} (x,y)$. \end{property}

To illustrate the computational limitations of this approach when $(SDP)$ is solved by standard semidefinite programming solvers, we refer to the experiments presented in Section~\ref{sec:kcluster} where we compare methods \texttt{MIQCR} and \texttt{QCR} on 45 $k$-cluster instances of size $n=80$. Recall that the $k$-cluster problem can be formulated by an equality constrained binary quadratic program. 

While \texttt{QCR}~\cite{BEP09} is devoted to solve binary quadratic programs, \texttt{MIQCR}~\cite{BEL12} was devised as a solution method for quadratic programs having general integer variables and continuous variables. In the end, the two methods are based on the same ideas of quadratic convex reformulation. \texttt{QCR} is in fact a special case of method \texttt{MIQCR} where all parameters $\beta_{ij}$ are fixed to 0.  Hence,  Phase~1 of \texttt{QCR} amounts to solve $(SDP)$ without Constraints~(\ref{ij1})-(\ref{ij4}) and Constraints~(\ref{i1}) replaced by $X_{ii}=x_i$. Moreover, in \texttt{QCR} no additional variables $y$ are necessary to get an equivalent formulation. Indeed, it can be seen from set $S_{xyzt}$ that variables $y_{ij}$, $i \neq j$, are created only if  $\beta_{ij} \neq 0$, and since $x_i \in \{0,1 \}$, we have $y_{ii} = x^2_i = x_i$. 

In \texttt{MIQCR} the equivalent formulation is thus of larger size than the initial problem. Naturally, the equivalent formulation obtained by \texttt{MIQCR} leads to a better bound than the one obtained by \texttt{QCR}~\cite{Lambert09}. Hence, \texttt{MIQCR} relies on a tighter semidefinite relaxation but with a larger number of additional constraints. Summarizing, in binary quadratic programming methods \texttt{QCR} and \texttt{MIQCR} differ by the compromise between tightness of the used bounds versus size of the semidefinite relaxation and of the reformulated problem. These sizes are clearly in relation with the cpu time needed for solving both, the semidefinite relaxation and the reformulated problem.

\section{\texttt{MIQCR-CB} - a method for large-scale problems}
\label{sec:solving}
In this section we propose a non-standard algorithm to compute the reformulated problem provided by method \texttt{MIQCR}. This algorithm is much faster than a standard semidefinite programming solver. Moreover, it allows us to start Phase~2 with a tight bound and a reformulated problem of reduced size, and thus to handle large instances.

As already explained, the bottleneck of method \texttt{MIQCR} is solving $(SDP)$.
The most prominent methods for solving semidefinite problems are interior-point methods, e.g.~\cite{HeReVaWo:96}. These methods are well-studied and several implementations exist, e.g., CSDP~\cite{csdp}, SeDuMi~\cite{sedumi}, SDPA~\cite{sdpa}. The computational effort depends on the order of the matrix and on the number of constraints. For instances with  matrix size larger than 1000, or with more than 10\;000 constraints, the semidefinite problems become intractable for interior-point methods.
In Section~4.1 we demonstrate the weakness of interior-point methods applied to our problem. It turns out that for $k$-cluster problems of size $n=80$ it is already not practical. 

A variety of alternative methods for solving semidefinite problems has been developed in the last decades. Many of these are based on augmented Lagrangian methods. In~\cite{BurMon03,BurMon05} an augmented Lagrangian algorithm is presented where the constraint $X\succeq 0$ is replaced by $X=RR^\top$, $R$ being a matrix of low rank.
Augmented Lagrangian algorithms using projection techniques are proposed in~\cite{JarRen08,MPRW09,ZST10}.

Another class for solving semidefinite programs are bundle methods.
In the spectral bundle method~\cite{specbdl} the semidefinite problem is reformulated as an eigenvalue optimization problem, which is then solved by a subgradient method. An implementation of this algorithm is SBMethod~\cite{sb}, or more general, the callable Conic Bundle library~\cite{cb}.

We choose to use a bundle method to obtain a reasonable solution within short time. Following the idea of~\cite{fischeretal}, we design subgradient algorithm within a Lagrangian duality framework.

\subsection{A static bundle method for solving $(SDP)$}
\label{subsec:static}
Let us consider a partial Lagrangian dual of $(SDP)$ where we dualize the linearization constraints, i.e. constraints~(\ref{ij1})--(\ref{ij4}). We rewrite $(SDP)$ as $(SDP_{T})$ using the following notation.
\begin{numcases}{(SDP_{T})}
\max  f(X,x) = \displaystyle{\sum_{i=1}^n} \displaystyle{\sum_{j=i+1}^n}  q_{ij} X_{ij} + \displaystyle{\sum_{i=1}^n}c_ix_i \nonumber \\
\mbox{s.t.} \nonumber \\
\quad (X,x) \in S  \nonumber\\
\quad h^t_{ij} (X,x) \leq 0, & $(i,j,t) \in T$ \nonumber
\end{numcases}

\noindent where $S = \{ (X,x)\colon (X,x) \mbox{ satisfies } (\ref{equality}), (\ref{c1})-(\ref{i3}), (\ref{sdp_1}), (\ref{sdp_2}) \}$ and  $T = \{(i,j,t): 1 \le i<j \le n, t=1, \ldots, 4 \}$, and for all $(i,j,t) \in T$:
\begin{numcases}{h^t_{ij}(X,x)=}
X_{ij} - u_jx_i  & $t=1$ \nonumber \\ 
X_{ij} - u_ix_j  & $t=2$ \nonumber \\
- X_{ij} + u_jx_i + u_ix_j -u_iu_j  & $t=3$ \nonumber\\
- X_{ij}      & $t=4$ \nonumber
\end{numcases}
With each constraint $h^t_{ij}(X,x) \leq 0$ of $(SDP_T)$ we associate a non-negative Lagrange multiplier $\beta^t_{ij}$. We now consider the partial Lagrangian 
$$\mathcal{L}_T(X,x,\beta) = \displaystyle{\sum_{i=1}^n} \displaystyle{\sum_{j=i+1}^n}  q_{ij} X_{ij} + \displaystyle{\sum_{i=1}^n}c_ix_i - \displaystyle{\sum_{(i,j,t) \in T}} \beta^t_{ij} h^t_{ij}(X,x)$$
and we obtain the dual functional
$$g_T(\beta) = \displaystyle{\max_{\scriptsize{(X,x) \in S}}} \quad \mathcal{L}_T(X,x,\beta).$$
By minimizing this dual functional we obtain the partial Lagrangian dual problem $(LD_T)$ associated with $(SDP_T)$, 
$$(LD_T) \left\{ \begin{array}{ll} 
\min & g_T(\beta)  \\
\mbox{s.t.} & \beta^t_{ij} \geq 0,\; (i,j,t) \in T. \\
\end{array} \right.$$
 Our aim is to solve $(LD_T)$ using the bundle method. The outline of the algorithm is the following. For a given $\bar{\beta} \geq 0$, we evaluate $g_T(\bar{\beta})$ and determine the associate primal solution $(\bar{X},\bar{x})$, such that $g_T(\bar{\beta}) = \mathcal{L}_T(\bar{X},\bar{x},\bar{\beta})$. We call a pair $(\bar{\beta}, (\bar{X},\bar{x}))$ a matching pair for $g_T$. Evaluating function $g_T$ for given $\bar{\beta}$ amounts to maximize a linear function in $(X,x)$ over the set $S$. This is an SDP that has much less constraints than $(SDP_{T})$ and can be solved efficiently by interior-point methods. From the solution $(\bar{X},\bar{x})$, we compute a subgradient $h_{ij}^t(\bar{X},\bar{x}) \in \partial g_T(\bar{\beta})$. The bundle method is an iterative algorithm that maintains at each iteration a ``best'' approximation $\hat{\beta}$ and a sequence $\mathcal{X}= ((\bar{X}_1,\bar{x}_1),(\bar{X}_2,\bar{x}_2), \ldots , (\bar{X}_k,\bar{x}_k))$ where $(\hat{\beta}, (\bar{X}_i,\bar{x}_i))$ is a matching pair. Then, from the sequence $\mathcal{X}$, the best approximation, $\hat{\beta}$, and the new subgradient, the bundle method computes a new value $\hat{\beta}$ that will be used at the next iteration. A detailed description of the method is available in~\cite{fischeretal}.

\subsection{A dynamic bundle method for solving $(SDP_T)$}
\label{subsec:dynamic}
In order to preserve efficiency we adopt another idea from~\cite{fischeretal}. The number of elements in $T$ is $4{n \choose 2}$. However, we are interested only in the subset of $T$ for which the constraints $h_{ij}^t(X,x) \leq 0$ are likely to be active at the optimum. This set is not known in advance, however, in the course of the algorithm we dynamically add and remove elements in order to identify ``important'' constraints. Here, we consider $\mathcal{T}\subseteq T$ and work with the function
$$g_\mathcal{T}(\beta) = \displaystyle{\max_{\scriptsize{(X,x) \in S}}} \quad \mathcal{L}_\mathcal{T}(X,x,\beta).$$
Initially we set $\mathcal{T}=\emptyset$ and after a first function evaluation we separate violated inequalities and add the elements to set $\mathcal{T}$ accordingly. We keep on updating this set in course of the bundle iterations by removing elements with associated multiplier close to zero and separate newly violated constraints. In this way we obtain a ``good'' set of constraints.

Convergence for dynamic bundle methods has been analyzed in detail in~\cite{BeSa:09}, giving a positive answer for convergence properties in a rather general setting.

\subsection{A parameterized dual heuristic for solving $(SDP_T)$ }
\label{subsec:heur}
The computation of the ``nearly'' optimal $(\alpha^*,\lambda^*,\beta^*)$ with the dynamic bundle method still can require much computational time. An idea for reducing this computational time is to consider a relaxation of $(SDP_T)$. Indeed, as observed in Section~2, any feasible dual solution to $(SDP_T)$ allows us to build a convex equivalent formulation to $(QP)$. A possible way to get such a solution is to drop some constraints from~(\ref{ij1})--(\ref{ij4}) of  $(SDP_T)$ and compute a dual ``nearly'' optimal solution $(\bar{\alpha},\bar{\lambda},\bar{\beta})$ of the reduced problem.  Then, a feasible dual solution to $(SDP_T)$ can be obtained by completing $(\bar{\alpha},\bar{\lambda},\bar{\beta})$  with zeros for those dual variables corresponding to the dropped constraints. To carry out this idea, we consider a parameter $p$ that is an upper bound on the cardinality of $\mathcal{T}$ ($\lvert \mathcal{T} \rvert \leq p$). In other words, $p$ is the maximum number of constraints considered in the reduced problem. Finally, the proposed dual heuristic has two extreme cases: 
\begin{itemize}
\item if $p= 4{n \choose 2}$, we solve $(SDP_T)$ and get the associated dual solution as in Section~3.2.
\item if $p=0$,  we make a single iteration: we get the optimal solution of the reduced problem obtained from $(SDP_T)$ where we drop all constraints~(\ref{ij1})--(\ref{ij4}) (For binary quadratic programming this amounts to method \texttt{QCR}).
\end{itemize}
 
\bigskip
 We call this procedure ComputeBeta($g_T$,$p$) and sketch it in Algorithm~\ref{algo1}. The algorithm returns a solution $\beta^*$ having at most $p$  positive components. Thus, the number of variables $y_{ij}$ of problem $(QP_{\alpha,\lambda,\beta})$ is also at most $p$ only and Phase~2 of \texttt{MIQCR-CB} can be solved much faster than Phase~2 of \texttt{MIQCR}. Finally, this parameter $p$ controls the size, and in a sense the tightness, of the semidefinite relaxation used for computing the equivalent formulation of method \texttt{MIQCR-CB}.

\begin{algorithm}
\begin{center}
\caption{ComputeBeta($g_T$,$p$)}\label{algo1}
\begin{algorithmic}[1]
\REQUIRE
\STATE k=0  \COMMENT{counter on the number of iterations}
\STATE $\mathcal{T}_k = \emptyset $ \COMMENT{Current set of dualized constraints : we start with $0$ for all $\beta$}
\STATE Solve $g_{\mathcal{T}_k}(0)$. Let $(\bar{X},\bar{x})$ be the obtained solution. \COMMENT {We start by solving $(SDP_{\mathcal{T}})$ without any Constraints~(\ref{ij1})--(\ref{ij4}) and determine the $p$ most violated constraints}
\STATE $\mathcal{T}_k = \{$ sub-set of $T$ corresponding to the at most $p$ violated constraints of $(SDP_T)$ at point $(\bar{X},\bar{x}) \}$.
\STATE Compute a subgradient $h_{ij}^t(\bar{X},\bar{x})$ of $g_T(\beta)$ for $\beta=0$.
\WHILE {predicted progress of the next step is sufficient}
\STATE Update $\beta$ value: $\beta = \hat{\beta}$ \COMMENT {using the Conic Bundle algorithm~\cite{cb}}
\STATE Solve $g_{\mathcal{T}_k}(\hat{\beta})$. Let $(\bar{X},\bar{x})$ be the obtained solution. \COMMENT{where the objective function of $g_{\mathcal{T}_k}(\hat{\beta})$ is the function obtained by dualizing each constraint $(i,j,t) \in \mathcal{T}_k$}
\STATE Compute a subgradient $h_{ij}^t(\bar{X},\bar{x})$ of $g_T(\beta)$ for $\beta=\hat{\beta}$.
\STATE Compute $\mathcal{T}_{k+1}$ : drop from $\mathcal{T}_k$ constraints that are no longer violated and add new most violated constraints such that $\lvert \mathcal{T}_{k+1} \rvert \leq p$ 
\STATE k++;
\ENDWHILE
\STATE Complete the solution $\beta^*$ by zeros for constraints that do not belong to $\mathcal{T}_k$.
\end{algorithmic}
\end{center}

\begin{remark}\label{remark-gbeta}
In Steps~3 and~8, $g_{\mathcal{T}_k}(\hat{\beta})$ is computed by CSDP~\cite{csdp}. 
\end{remark}

\end{algorithm}

\section{Computational results}
\label{sec:comp}

  In this section, we present computational results for our method \texttt{MIQCR-CB}. We first evaluate our algorithm on binary quadratic programming instances of the $k$-cluster problem. For this, we start with a detailed comparison of methods \texttt{QCR}, \texttt{MIQCR}, and  \texttt{MIQCR-CB} for instances of size $n=80$. For these instances, we also study the behavior of \texttt{MIQCR-CB}  when varying parameter $p$. As instances of size $n=80$ are not practical for \texttt{CSDP}, we make experiments using instances of smaller size ($n=40$) to compare \texttt{MIQCR-CB} with the interior-point solver \texttt{CSDP} turned into a heuristic. These experiments illustrate that for Phase~1 the bundle algorithm is faster. Finally, we compare our method with \texttt{BiqCrunch}~\cite{KMR13} for larger instances of sizes $n\in \{100,120,140,160\}$. Furthermore, we evaluate \texttt{MIQCR-CB} on instances with general integer variables. For theses instances, we compare our algorithm with the solvers \texttt{Couenne}~\cite{couenne} and \texttt{cplex~12.6}~\cite{cplex126} as the scope of \texttt{BiqCrunch} is binary quadratic programming only. Note that our method can also handle general mixed-integer problems, some computational results can be found in~\cite{BEL15}.\\

\noindent \textbf{Experimental environment:} \\

We implemented algorithm \texttt{MIQCR-CB} in~C. We use the \textit{Conic Bundle} callable C-library of Christoph Helmberg~\cite{cb} to implement Algorithm~1 and the SDP solver CSDP of Brian Borchers~\cite{csdp} for the function evaluation. Methods \texttt{MIQCR} and \texttt{QCR} are also available as C implementations.
For solving $(SDP)$ of method \texttt{MIQCR} we use the solver SBMethod~\cite{sb}, as the solver CSDP~\cite{csdp} was not able to handle $(SDP)$ (allocation storage error), and the solver CSDP~\cite{csdp} is used for solving the semidefinite programs of \texttt{QCR}. The C-interface of \texttt{Cplex12.5}~\cite{cplex125} serves for solving the quadratic programs.

Experiments for all methods (\texttt{MIQCR-CB}, \texttt{MIQCR}, \texttt{QCR} and \texttt{BiqCrunch}) were carried out on a laptop with an Intel quad-core $i7$ processor of 1.73 GHz and $6$ GB of RAM using a Linux operating system. \\

\subsection{Computational results for the $k$-cluster problem}
\label{sec:kcluster}

Given a graph $G$ of $n$ vertices and a number $k \in \{3, \ldots, n-2\}$, the $k$-cluster problem consists in finding a subset of $k$ vertices of $G$ such that the induced subgraph is as dense as possible. This problem or its weighted version has many applications and is classical in combinatorial optimization. It is also known as the ``heaviest $k$-subgraph problem'', the ``$k$-dispersion problem''~\cite{Pisinger06}, the ``$k$-defense-sum problem''~\cite{Krarup02}, the ``densest $k$-subgraph problem'' and the ``$k$-subgraph problem''. It can be formulated by the following binary quadratic program.
\begin{numcases}{(KC)}
\max  f(x)= \displaystyle{\sum_{i=1}^n} \displaystyle{\sum_{j=i+1}^n}  \delta _{ij} x_i x_j \nonumber \\
\mbox{s.t.} \nonumber \\
\quad \displaystyle{\sum_{i=1}^n} x_i= k  \nonumber\\
 \quad x \in \{0,1\}^n  \nonumber
\end{numcases}
where $\delta _{ij}= 1$ if and only if an edge links vertices $i$ and $j$ and $x_i$ is the binary variable indicating whether vertex $i$ is selected in the subset. $(KC)$ is known to be $\mathcal{NP}$-hard even for bipartite graphs~\cite{Corneil84}. Many approximation results are known for $(KC)$ ~\cite{Asahiro96,Hassin97,Kortsarz93,Srivastav98}. Concerning the solution algorithms of $(KC)$, classical approaches based on linearization techniques are able to solve medium size instances with up to 80 variables~\cite{Bil05,Erkut90,Pisinger06}. A few methods are able to solve $(KC)$ to optimality for large size instances (when $n>80$). The most efficient exact solution methods are based on nonlinear approaches, such as convex quadratic programming~\cite{BEP09} or semidefinite programming~\cite{HanYe02,JagSri05,MalRou11,Roupin04}.

We compare experimentally our new algorithm \texttt{MIQCR-CB} with three approaches: the original \texttt{MIQCR} and \texttt{QCR} approaches, and the method \texttt{BiqCrunch} of Krislock, Malick, and Roupin~\cite{KMR13}. This latter approach uses the branch-and-bound solver Bob~\cite{bob} in monothreading together with semidefinite programming bounds~\cite{Mal07} to solve $(KC)$ to optimality. At each node, a dual bound is computed solving a semidefinite relaxation of $(KC)$. The semidefinite relaxation turns out to be tighter than the one used in \texttt{MIQCR} because it integrates the family of triangular inequalities that are not used inside \texttt{MIQCR}. To compute the dual bound, the SDP program is first formulated as an equivalent program with a spherical constraint (a constraint on the norm of matrix $X$). This spherical constraint is then dualized within a Lagrangian framework, and the dual problem is viewed and solved as a particular least-squares semidefinite program. The resulting bound is very close to the optimal value of the semidefinite relaxation and can be computed very fast.

We report numerical results on $225$ instances of $(KC)$ with up to $160$ vertices. We use the 90 instances of sizes $n=80$ and $100$ introduced in~\cite{Bil05}, extended in~\cite{BEP05} and~\cite{Plateau06}, and also used in~\cite{KMR13}. Additionally we consider the set of 135 instances of sizes $n\in \{120,140,160\}$ used in~\cite{KMR13}. All these instances were generated as follows. For a given number of vertices $n$ and a density $d$ an unweighted graph is randomly generated. The parameter $k$ is then set to $\frac{n}{4}$, $\frac{n}{2}$, and $\frac{3n}{4}$. All instances are available online~\cite{lib_kcluster}.\\

\noindent \textbf{Parameters:}
\begin{itemize}
\item Phase~1: Parameter \texttt{termeps} of SBMethod~\cite{sb} is set to $10^{-4}$. Parameters \texttt{axtol, aytol} of CSDP~\cite{csdp} are set to $10^{-4}$. The precision of the \textit{Conic Bundle}~\cite{cb} is set to $10^{-4}$. For method \texttt{MIQCR-CB}, we allow to dualize all constraints, i.e., $p=4{n \choose 2}$.

\item Phase~2: The tolerance for parameter $\beta_{ij}$ to be considered as non-zero is $10^{-4}$. For \texttt{Cplex12.5}~\cite{cplex125} (used in \texttt{QCR},  \texttt{MIQCR}, and  \texttt{MIQCR-CB}), the relative mipgap is $10^{-6}$ and the absolute gap is $0.99$. The parameter \texttt{objdiff} is set to $0.999$, and the parameter \texttt{varsel} to $4$. The time limit is set to 3~hours. We use the multi-threading version of \texttt{Cplex12.5} with up to 8 threads. 
\end{itemize}

\noindent \textbf{Legend of Tables~\ref{tab:kcluster80} and~\ref{tab:bc140}--\ref{tab:ieqp}}
\begin{itemize}
\item Each line is an average over 5 instances and we consider 45 instances for each size 80, 100, 120, 140, 160;
\item $n$ indicates the size of the graph;
\item  $k$ is the size of the subgraph or cluster;
\item  $d$ is the density of the graph;
\item \textit{Gap} is the relative gap in percentage between the optimal solution value $v$ and the value of the continuous relaxation $c$ at the root node of the branch-and-bound tree ($Gap\, = \, \frac{\lvert c - v \rvert}{v}*100$);
\item \textit{P1} is the cpu time in seconds for solving Phase~1;
\item \textit{P2} is the cpu time in seconds for solving Phase~2 by \texttt{Cplex12.5}~\cite{cplex125}; \texttt{(i)}~means that only \texttt{i} instances were solved within the time limit, averages are taken only over instances solved within the time limit;
\item \textit{Tt} is the total time in seconds (the sum of \textit{P1} and \textit{P2} for \texttt{MIQCR} and \texttt{MIQCR-CB}). If the optimum is not found within this time, we present the final
  gap $(g \%)$, $g = \displaystyle \left| \frac {Opt - b }{ Opt} \right| * 100$ where $b$ is the best bound obtained within the time limit;
\item \textit{Min} and \textit{Max} are the minimum and maximum total time, respectively, within the five instances of the same characteristics. 
\item  \textit{Nodes} is the number of nodes explored.\\
\end{itemize}

%\begin{table}
\begin{table}
\centering
\begin{scriptsize}
\begin{tabular}{c|c|c||c|c|c|c|c||c|c|c|c|c||c|c|c|c|c}
\multicolumn{3}{c||}{}&  \multicolumn{5}{|c||}{\texttt{MIQCR}}  & \multicolumn{5}{|c||}{\texttt{QCR $(p= 0)$}} & \multicolumn{5}{|c}{\texttt{MIQCR-CB ($p= 4{n \choose 2}$)}} \\  
 \textit{n}&\textit{d $(\%)$}&\textit{k} & \textit{Gap} &  \textit{P1} &  \textit{P2}&  \textit{Tt}  & \textit{Nodes} & \textit{Gap} &  \textit{P1} &  \textit{P2} &  \textit{Tt}   & \textit{Nodes}  & \textit{Gap} &  \textit{P1} &   \textit{P2}&  \textit{Tt} & \textit{Nodes} \\\hline
80&25&20&3.39&	1434&	129&	1563 & 2371&  9.2&1&12&13&77084 &3.00&6&4&10   &940\\
80&25&40&1.00&482&	83&	565&	924&2.68&1&3&4&12787&0.72&6&3&9    &204\\
80&25&60&0.30&199&	39&	238&	182 &0.87&1&2&3&1137&0.07&5&3&8   &0\\\hline
80&50&20&2.34&	981&	173&	1154&	3281&8.03&1&22&23&231260&2.02&6&5&11  &2013\\
80&50&40&0.76&	373&	143&	516&	1809&1.81&1&3&4&14049&0.49&6&3&9   &426 \\
80&50&60&  0.29	&178&	215&	393&	1748&0.6&1&2&3&3249&0.08&5&3&8   &8\\\hline
80&75&20&1.49&1273&	188&	1461& 3912&6.47&1&66&67&792944&1.37&5&6&11  &2081 \\
80&75&40&0.59&411& 988	&1399	&19621&1.35&1&12&13&113410&0.49&6&6&12  &5353\\
80&75&60&0.20&	220&	132&	352&	917&0.42&1&2&3&3247&0.04&5&3&8  &0\\\hline
 \multicolumn{3}{c||}{\textbf{Mean} } & \textbf{1.15} & \textbf{617} & \textbf{232} & \textbf{849}&  \textbf{3863} & \textbf{3.49}&\textbf{1}&\textbf{14}&\textbf{15}&\textbf{138796}&\textbf{0.92}&\textbf{6}&\textbf{4}&\textbf{10}    &\textbf{1225}

\end{tabular}
\end{scriptsize}
\caption {Average computational results of \texttt{MIQCR}, \texttt{QCR} and  \texttt{MIQCR-CB} for 45 $k$-cluster instances with $n=80$.}  
\label{tab:kcluster80}
\end{table}
To illustrate how the three methods compare with respect to bound, cpu time and nodes, we show in Table~\ref{tab:kcluster80} numerical results for 45~instances of size $n=80$ solved using \texttt{MIQCR}, \texttt{QCR} and  \texttt{MIQCR-CB}, where each line corresponds to average values over $5$ instances, the bottom line averages the values over all instances. \\

\noindent \textbf{Comparison of \texttt{MIQCR} and \texttt{QCR}\\}

As expected, we observe that the gap obtained by \texttt{MIQCR} is always smaller than the gap obtained by \texttt{QCR}. On average the gap is only one third of the gap obtained by \texttt{QCR}. In spite of this improvement of the gap, the cpu time of  \texttt{MIQCR} is $57$ times larger than the cpu time of \texttt{QCR}. This is caused by the addition of the variables $y$ and the corresponding linearization constraints in the reformulated problem. Indeed, in \texttt{QCR} the reformulation has $n$ variables and $m$ constraints, while in \texttt{MIQCR} the reformulation has in the worst case (i.e. all $\beta_{ij} \neq 0$) $n + n\lvert E\rvert + \lvert E\rvert$ variables and $5n^2 + 4n\lvert E\rvert +4n$ constraints. Observe, however, that as a consequence of the smaller root gap, about $36$ times less nodes are explored by \texttt{Cplex12.5} in \texttt{MIQCR}.\\

\noindent Concerning the comparison of \texttt{MIQCR} and \texttt{QCR}, we summarize:

\begin{enumerate}[i)]

\item The reformulated problem for method~\texttt{QCR} is faster to compute, but the theoretical continuous relaxation bound provided by its equivalent formulation is weaker. 
As a consequence, the solution computed by the associated branch-and-bound algorithm in Phase~2 of \texttt{QCR} is hard to obtain. 

\item Conversely, the reformulated problem provided by~\texttt{MIQCR} leads to a tight theoretical continuous relaxation bound, but this equivalent formulation computed in Phase~1 is very hard to solve. 

\item For both methods, the solution of {\em large-scale} $k$-cluster problems is intractable. 
\end{enumerate}

\noindent \textbf{Comparison of \texttt{MIQCR-CB} and \texttt{MIQCR}}\\

We would like to demonstrate the improvement of \texttt{MIQCR-CB} over \texttt{MIQCR}. When comparing these two latter methods, we can observe the following.
\begin{enumerate}[i)]

\item \textit{The computation time of Phase~1 is significantly reduced.} Indeed, the average computational time over all the instances reduces from $617$ seconds to $6$ seconds. 

\item \textit{The SDP bound is tightened.} For these instances, the gap obtained by \texttt{MIQCR-CB} is smaller compared to \texttt{MIQCR}. This might sound strange because the SDP problem considered is the same for both methods. The difference results from our algorithm used to solve Phase~1 that is more accurate than SB for these instances.

\item \textit{The computation time of Phase~2 is reduced.} This time is divided by a factor $58$ on average for \texttt{MIQCR-CB}. This is mainly due to the smaller number of variables in the reformulated problem. Indeed, we have a significant number of $\beta_{ij}$ values that are $0$, and thus less variables $y_{ij}$ with their associated constraints are considered in the reformulated problem. As an illustration of this phenomenon, for one of the instances with $n=80$, $d=25\%$, and $k=20$, we observe $2488$ non-zero $\beta_{ij}$ in \texttt{MIQCR}, and only $1890$ non-zero $\beta_{ij}$ in \texttt{MIQCR-CB}.\\
\end{enumerate}

\noindent \textbf{Comparison of \texttt{MIQCR-CB} and \texttt{QCR}}\\

 We make the following observations.

\begin{enumerate}[i)]

\item \textit{P1 is smaller for} \texttt{QCR} \textit{than for} \texttt{MIQCR-CB}.
For \texttt{QCR}, \textit{P1} is always smaller than one second. In \texttt{MIQCR-CB} \textit{P1} varies from 5 to 6~seconds. This is the consequence of adding Constraints~(\ref{ij1})-(\ref{ij4}) which number is potentially $\Theta(|I|)$. 

\item \textit{The gap in }\texttt{MIQCR-CB} \textit{is significantly smaller than the gap in} \texttt{QCR}.
Compared to \texttt{QCR}, the gap in \texttt{MIQCR-CB} is approximatively divided by 1.5 for all the instances. This is also the consequence of using a stronger semidefinite relaxation. 

\item \textit{On average} \texttt{MIQCR-CB} \textit{is faster than} \texttt{QCR}.
The average total cpu time of \texttt{MIQCR-CB} is divided by a factor $4$ compared to \texttt{QCR}. Note that for both methods the $k$-cluster problem gets harder with smaller values of $k$. Especially for small values of $k$, \texttt{MIQCR-CB} is significantly faster than \texttt{QCR}, although \texttt{QCR} is faster for these medium size instances for half cases.

\end{enumerate}
All the instances considered of size $n=80$ are solved to optimality by \texttt{MIQCR-CB} within 20 seconds. We present in the next section results for larger problems, to determine the limitations of method  \texttt{MIQCR-CB} and compare it with \texttt{BiqCrunch}. \\

%%%%%%%%%%%%%%%%%%%%%%%%%%%%%%%%%%%%%%%%%%%%%%%%%%%%%%%%%%%%%%%%%%%%%%%%%%%%%%%%%%%%%%%%%%%%%%%%%%%%%%%%%%%%%%%%%%%%%%%%%%%%%%%%
%%%%%%%%%%%%%%%%%%%%%%%%%%%%%%%%%%%%%%%%%%%%%%%%%%% n = 80 %%%%%%%%%%%%%%%%%%%%%%%%%%%%%%%%%%%%%%%%%%%%%%%%%%%%%%%%%%%%%%%%%%%%%
%%%%%%%%%%%%%%%%%%%%%%%%%%%%%%%%%%%%%%%%%%%%%%%%%%%%%%%%%%%%%%%%%%%%%%%%%%%%%%%%%%%%%%%%%%%%%%%%%%%%%%%%%%%%%%%%%%%%%%%%%%%%%%%%

\noindent \textbf{Computational study of the influence of parameter $p$\\}

 In Section~3.3, we introduced a dual heuristic version of our bundle algorithm used for solving Phase~1 that is parameterized by $p$. This parameter controls the size, and the tightness, of the semidefinite relaxation used for computing the equivalent formulation of method \texttt{MIQCR-CB}. To evaluate the influence of $p$ on method  \texttt{MIQCR-CB}, we run our method for different values of $p$. Denote by $p_0= 4{n \choose 2} =\lvert T \rvert$ the initial number of inequalities (\ref{ij1})--(\ref{ij4}), we run our method for $p= \delta p_0$, with $\delta = 1, 0.5, 0.2,0.1,0.05$ and $0.01$. We report in Table~\ref{tab:p} the corresponding values of $\delta$ and of $p$ for these instances of size $80$. We can see in Figure~1 the results obtained for the initial gap, the solution time of Phase~1, and the total solution time. We observe that if the number of active constraints at the optimum is strictly smaller than $p$, the optimal solution of $(SDP)$ is obtained faster when $\delta$ is set to 1. Indeed, see for instance when $\delta=0.5$, the solution time is larger than for $\delta=1$, and the initial gaps are the same for both parameters. Otherwise, the smaller the $p$, the faster the solution time of $(SDP)$ is computed. These results also illustrate the advantage of method   \texttt{MIQCR-CB} over method  \texttt{QCR} (when $p=0$), as for any considered values of $p$, \texttt{MIQCR-CB} is always faster.

\begin{table}
\centering

\begin{tabular}{|c|c|}\hline
 $\delta$ & $p$ \\ \hline
 $ 1$ & 3160  \\ \hline
 $0.5$& 1580 \\ \hline
 $0.2$  & 632  \\ \hline
 $0.1$ & 316  \\ \hline
 $0.05$ &158  \\ \hline
 $0.01$ & 31  \\ \hline

\end{tabular}
\caption {Number $p$ of considered constraints with different values of $\delta$}  
\label{tab:p}
\end{table}

\begin{figure}
\label{comp_p}
\begin{center}
 \includegraphics[width=10cm]{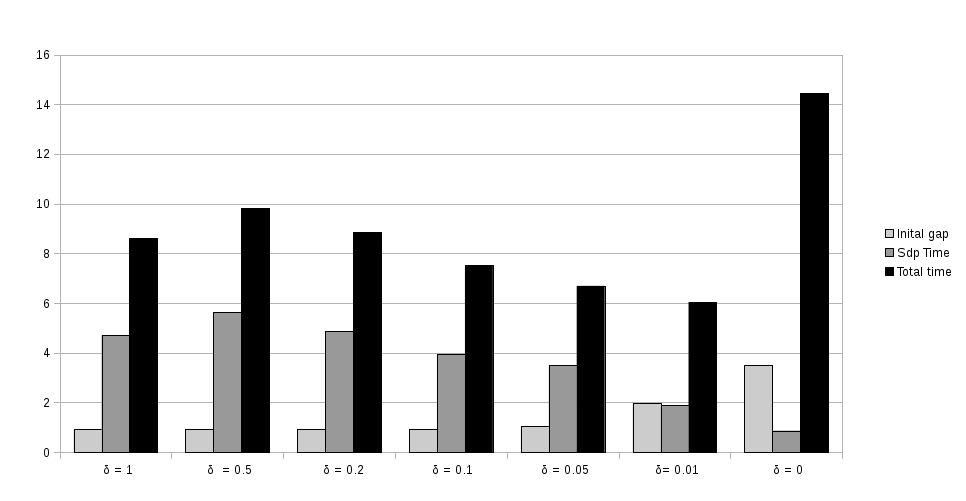}
\end{center}
\caption{$k$-cluster instances of size  $n=80$ : influence of parameter $p$}
\end{figure}
\noindent \textbf{Experimental comparison of \texttt{MIQCR-CB} (with $p =  4{n \choose 2}$) and \texttt{MIQCR} where Phase~1 is solved with \texttt{CSDP} turned into a heuristic\\}

In this paper, we want to accelerate Phase~1 of our algorithm, and furthermore get an equivalent formulation with a reasonable size. A way to do this using interior-point algorithm for solving $(SDP)$ is to stop the solver \texttt{CSDP} with a much smaller precision and to consider the $\beta$ equals to zero with a much smaller precision. As instances of size $n=80$ are not practical for \texttt{CSDP}, we make experiments using instances of smaller size ($n=40$). We use 3 configurations for solving Phase~1 of our method:
\begin{itemize}
\item \textit{Configuration 1: } We use \texttt{CSDP} turned into a heuristic: parameters \texttt{axtol} and \texttt{aytol}  are set to $10^{-4}$ and parameter \texttt{objtol} is set to $10^{-1}$, and the tolerance for parameters $\beta_{ij}$ to be considered as non-zero to $10^{-4}$. 
\item \textit{Configuration 2: } We use \texttt{CSDP} turned into a heuristic: parameters \texttt{axtol}, \texttt{aytol}  are set to $10^{-4}$ and \texttt{objtol} is set to $10^{-1}$, and the tolerance for parameters $\beta_{ij}$ to be considered as non-zero to $10^{-2}$. 
\item \textit{Configuration 3: } We use Algorithm~1: the stopping criteria is set to $10^{-4}$, and the tolerance for parameters $\beta_{ij}$ to be considered as non-zero to $10^{-4}$. We set $p$ to $4{n \choose 2}$.
\end{itemize}

We report in Table~\ref{tab:csdp} the average number of variables $y_{ij}$ created in Phase~2 for each considered configuration of Phase~1. We recall that for these instances the maximum number of variables $y_{ij}$ created is 820.  

\begin{table}
\centering

\begin{tabular}{|c|c|}\hline
Config. for Phase~1 & nb $y_{ij}$ in Phase~2 \\ \hline 
Configuration 1 &  738\\ \hline
Configuration 2 & 554\\ \hline
Configuration 3 & 409\\ \hline
\end{tabular}

\caption {Average number of variables $y_{ij}$ considered in Phase~2 of our algorithm}  
\label{tab:csdp}
\end{table}

We report in Figure~2 the results obtained for the initial gap, the solution time of Phase~1, and the total solution time. In this graphic, the scale is logarithmic. We observe that the initial gap obtained with configuration~3 is about 22 times smaller than the initial gap obtained by configurations~1 and~2. Moreover, the solution time of Phase 1 is 68 times smaller for configuration~3, in comparison to configurations~1 and~2.  

\begin{figure}
\label{comp_csdp}
\begin{center}
 \includegraphics[width=6cm,height=6cm]{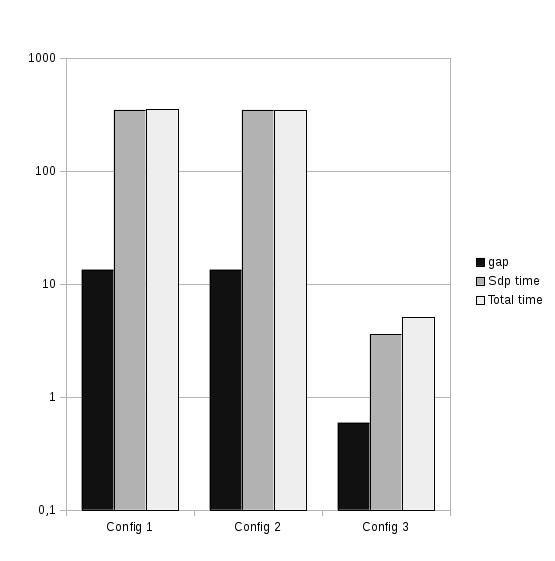}
\end{center}
\caption{$k$-cluster instances of size  $n=40$ : comparison with \texttt{CSDP} turned into a heuristic}
\end{figure}

Experiments to show the limitations of interior-points methods for solving problems with a huge number of constraints were already done in~\cite{fischeretal} and revealed a similar trend as ours. \\

\noindent \textbf{Comparison of \texttt{MIQCR-CB} and \texttt{BiqCrunch}\\}

In Table~\ref{tab:bc140} we report computational results on instances of larger size, namely $n \in \{100,120,140,160\}$. For these instances, we compare method \texttt{MIQCR-CB} with \texttt{BiqCrunch} run on the same computer.

 For instances of size $n=100$ and $n=120$, we observe that \texttt{MIQCR-CB} is as fast as \texttt{BiqCrunch}. Indeed, \texttt{MIQCR-CB} solves all $45$ instances of size $n=100$ ($n=120$ resp.) within $37.44$ ($370.89$ resp.) seconds on average, while \texttt{BiqCrunch} solves all the instances within $34.53$ ($338.07$ resp.) seconds. Notice that the initial gap is $1.6$ ($1.5$ resp.) times tighter for \texttt{BiqCrunch} than for \texttt{MIQCR-CB}. The average number of nodes is $1003$ ($1834$ resp.) times larger for \texttt{MIQCR-CB} for instances of size $n=100$ ($n=120$ resp.).

For the largest instances, namely those of size $n=140$ ($n=160$ resp.), \texttt{MIQCR-CB} is able to solve 38 (26 resp.) instances out of 45 within the time limit of 3~hours while \texttt{BiqCrunch} solves 43 (33 resp.) instances within 3 hours. The average computation time is 1973 (2081 resp.) seconds for \texttt{MIQCR-CB} and 1475 (1905 resp.) seconds for \texttt{BiqCrunch}. The initial gap is $1.4$ ($1.3$ resp.) times tighter for \texttt{BiqCrunch} than for \texttt{MIQCR-CB}. These largest instances represent the limit of both algorithms.

%%%%%%%%%%%%%%%%%%%%%%%%%%%%%%%%%%%%%%%%%%%%%%%%%%%%%%%%%%%%%%%%%%%%%%%%%%%%%%%%%%%%%%%%%%%%%%%%%%%%%%%%%%%%%%%%%%%%%%%%%%%%%%%%
%%%%%%%%%%%%%%%%%%%%%%%%%%%%%%%%%%%%%%%%%%%%%%%%%%% n = 100,120, 140, 160 %%%%%%%%%%%%%%%%%%%%%%%%%%%%%%%%%%%%%%%%%%%%%%%%%%%%%%
%%%%%%%%%%%%%%%%%%%%%%%%%%%%%%%%%%%%%%%%%%%%%%%%%%%%%%%%%%%%%%%%%%%%%%%%%%%%%%%%%%%%%%%%%%%%%%%%%%%%%%%%%%%%%%%%%%%%%%%%%%%%%%%%

\begin{table}
\centering
\begin{scriptsize}
\begin{tabular}{c|c|c||c|c|c|c|c|c|c||c|c|c|c|c}
\multicolumn{3}{c||}{}&  \multicolumn{7}{|c||}{\texttt{MIQCR-CB ($p= 4{n \choose 2}$)}} & \multicolumn{5}{|c}{\texttt{BiqCrunch}}  \\  
 \textit{n}&\textit{d (\%)}&\textit{k} & \textit{Gap}&  \textit{P1} &  \textit{P2} &  \textit{Tt} &  \textit{Min} &  \textit{Max}  & \textit{Nodes}& \textit{Gap} &  \textit{Tt}  &  \textit{Min} &  \textit{Max}& \textit{Nodes}\\\hline

100 & 25 &25 &3.30&14.4&25.2&39.6& 19 & 82 & 11964&1.78&71.3& 14 & 265 & 26\\
100 & 25 &50 &0.84&13.6&9.8&23.4& 17 & 44 &3407&0.28&9.3& 2&29&4\\
100 & 25 &75 &0.08&11.8&5.2&17&16&18&5&0.11&2.6&2&3&1\\
											
100 & 50 &25 &2.50&13.6&41.2&54.8&21&116&29116&1.76&59.4&19&140&33\\
100 & 50 &50 &0.90&17.2&87.8&105& 18& 176&77626&0.43&83.3& 3& 149&28\\
100 & 50 &75 &0.08&14.4&6.2&20.6& 16&28&91&0.06&3.9& 2& 11&1\\

100 & 75 &25 &1.37&13.6&21.4&35&23 &54& 12559&1.19&75.1& 44&116&43\\
100 & 75 &50 &0.35&13.8&10.2&24&19&31&3603&0.10&3.3& 2&7&1\\
100 & 75 &75 &0.04&11.8&5.8&17.6&17&19&31&0.03&2.6&2&5&1\\

%\multicolumn{3}{c||}{\textbf{Mean} } & \textbf{1.05} & \textbf{14} & \textbf{24} & \textbf{38} &&& \textbf{15378} & \textbf{0.64}& \textbf{35} &&& \textbf{15} \\
\hline

120 & 25&30&3.53&31.2&334.8&366&168&818&140208& 2.22&348.5&224&522&104\\
120 & 25&60&1.03&33.4&343&376.4&75&695&151256&0.53&285.3&52&582&69\\
120 & 25&90&0.10&29&9.8&38.8&34&48&37&0.03&9.2&4&17&2\\
120 & 50&30&2.53&29.2&444.8&474&208&867&232762&1.92&581.9&251&1161&186\\
120 & 50&60&0.81&32.6&601.2&633.8&310&1246&322930&0.41&360.4&249&446&85\\
120 & 50&90&0.14&34&56.2&90.2&50&160&2895&0.08&135.1&38&301&22\\
120 & 75&30&1.57&31.6&791.6&823.2&56&2922&576795&1.42&1124.9&74&3232&403\\
120 & 75&60&0.47&40.4&457.8&498.2&111&903&245271&0.21&193.4&66&392&40\\
120 & 75&90&0.04&29&8.4&37.4&32&46&79.4&0.02&4.2&3&6&1.00\\
%\multicolumn{3}{c||}{\textbf{Mean} } & \textbf{1.14} & \textbf{32} & \textbf{339}& \textbf{371} &&& \textbf{185803}&\textbf{0.76}&\textbf{338}&&& \textbf{101}\\
\hline

140 & 25&35&3.55&53.6&4458&4511.6&452&9940&1624107&2.41&1691.7&305&3493&349\\
140 & 25&70&1.07&61.2&2725.4&2786.6&221&7939&1124785&0.57&864.6&129&2456&143\\
140 & 25&105&0.09&46.2&15.6&61.8&51&75&134&0.04&19.2&5&25&3\\
140 & 50&35&2.61&54.3&3591.3(4) &3645.5 (4) &1262&7062&1501072&2.16&1958.9 (4)&806&3589&493\\
140 & 50&70&0.89&42&1871(1)&1913 (1)&1913&1913&1013112&0.54&4663.7&379&8763&800\\
140 & 50&105&0.08&50&43.8&93.8&54&151&1462&0.04&91.9&7&251&11\\
140 & 75&35&1.69&51.7&3748.3(3)&3800 (3)  &2000&7048&1924558&1.57&4289.9 (4)&2475&9031&1128\\
140 & 75&70&0.39&54.4&1667.8&1722.2&493&4843&700401&0.19&305.9&185&551&46\\
140 & 75&105&0.05&47.2&193.2&240.4&47&721&9100.6&0.03&49.2&4&140&6\\
%\multicolumn{3}{c||}{\textbf{Mean} } & \textbf{1.16} & \textbf{52} & \textbf{2035 (38)}& \textbf{2086 (38)} &&& \textbf{877637}&\textbf{0.84}&\textbf{1548 (43)}&&& \textbf{331}\\
\hline

160 & 25&40&3.27&78.4&3932.7 (3) &4009.3 (3) &1705&7986&769425&2.29&2874.2 (4)&829&7078&453\\
160 & 25&80&1.02&77.5&4049(2)&4126.5(2)&2117&6136&811861&0.60&2737.7&491&6662&357\\
160 & 25&120&0.16&84.4&256.4&340.8&71&739&8872&0.09&353.1&38&899&31\\
160 & 50&40&2.54&69&9845(1) &9914(1)&9914&9914&2768184&2.15&6041.1(2)&5103&6979&1037\\
160 & 50&80&0.66&83&4665(2)&4748(2)&2182&7314&982063&0.42&4058(4)&558&9369&536\\
160 & 50&120&0.08&78.8&220.2&299&86&264&4222&0.04&212.5&41&468&16\\
160 & 75&40&1.68&110&3483(1)&3593(1)&3593&3593&871354&1.55&4518(1)&4518&4518&657\\
160 & 75&80&0.46&88.5&2176(2)&2264.5(2)&379&4150&403871&0.29&643.0(2)&84&1202&40\\
160 & 75&120& 0.05& 66.8 &553 &619.8&100&1392&12461 &0.03&	151.4&42&281&	12\\

%\multicolumn{3}{c||}{\textbf{Mean} } & \textbf{1.10}& \textbf{78}& \textbf{ 3242 (26)} & \textbf{3324 (26)} & &&\textbf{736924}&\textbf{0.83}&\textbf{2549 (33)}&&& \textbf{349}\\

\end{tabular}
\end{scriptsize}
\centerline{\scriptsize{(i): i instances out of 5 were solved within the time limit. The reported values correspond to these instances.}}

\caption {Average computational results of \texttt{MIQCR-CB} and \texttt{BiqCrunch} for 180 instances of the $k$-cluster of size $n=100$, $120$, $140$ and $160$ (time limit: 10\;800 seconds).}
\label{tab:bc140}
\end{table}

\subsection{Computational results on the general integer case}
\label{sec:general}
In this section we extend our computational experiments to the class of (general) quadratic integer problems to compare the improvement of \texttt{MIQCR-CB} over \texttt{MIQCR} on this more general class of instances, and evaluate our method with the state-of-the-art solvers \texttt{Cplex12.6}~\cite{cplex126} and \texttt{Couenne}~\cite{couenne}. We briefly recall that the solver \texttt{Couenne}~\cite{couenne} uses linear relaxations within a spatial branch-and-bound algorithm~\cite{BLL09}, and that \texttt{Cplex12.6} uses a branch-and-bound algorithm based on convex relaxations~\cite{BliBon13}. For both solvers, numerous heuristics are incorporated into branch-and-bound algorithms in order to improve their performances.\\\\
\noindent \textbf{Parameters:} 
\begin{itemize}
\item Parameters \texttt{axtol, aytol} of CSDP~\cite{csdp} are set to $10^{-8}$. The precision of the \textit{Conic Bundle}~\cite{cb} is set to $10^{-8}$. For method \texttt{MIQCR-CB}, we allow to dualize all constraints (i.e. $p=4{n \choose 2}$).
\item Phase~2: A parameter $\beta_{ij}$ is considered as $0$ if $|\beta_{ij}|$ is below $10^{-6}$. For  \texttt{Cplex12.5}~\cite{cplex125}, the relative mipgap is $10^{-8}$ and the absolute gap is $0.99$. The parameter \texttt{objdiff} is set to $0.999$, and the parameter \texttt{varsel} to $4$. We used the multi-threading version of \texttt{Cplex12.5} and \texttt{Cplex12.6} with up to 8 threads. For solvers \texttt{Cplex 12.6} and \texttt{Couenne} we keep the default parameters.
\end{itemize}

\noindent \textbf{Experiments on the Equality Integer Quadratic Problem} \\

We consider the Equality Integer Quadratic Problem (EIQP) that consists in minimizing a quadratic function subject to one linear equality constraint:\\

\begin{numcases}{(EIQP)}
\min  f(x)=  \displaystyle{\sum_{i=1}^n} \displaystyle{\sum_{j=i}^n}  q_{ij} x_i x_j  + \sum_{i=1}^n c_ix_i  \nonumber \\
\mbox{s.t.} \nonumber \\
\quad \displaystyle{\sum_{i=1}^n}a_i x_i = b \nonumber\\
\quad 0 \leq x_i \leq u_i& $1\le i \le n$ \nonumber\\
\quad x_i \in \mathbb{N}& $1\le i \le n$ \nonumber
\end{numcases}

We use the instances introduced in~\cite{BEL12}, available online~\cite{eiqp}. We run experiments on two classes of problem, EIQP\textunderscore 1  and EIQP\textunderscore 2. For each class we generate instances with $20$, $30$, and $40$ variables where the coefficients are randomly generated as follows:
\begin{itemize}
\item The coefficients of $Q$ and $c$ are uniformly distributed integers from the interval $[-100,100]$. 
\item The $a_i$ coefficients are uniformly distributed integers from the interval $[1,50]$ for class~1 and from $[1,100]$ for class~2.
\item  $b= \mu\cdot \displaystyle{\sum_{i=1}^n} a_i$, where $\mu=15$ for class~1 and  $\mu=20$ for class~2.
\item $u_i=30$ for $1 \le i \le n$ for class~1, and $u_i=50$ for class~2.
\end{itemize}

For each class and each $n\in \{20,30,40\}$ we consider $5$ instances obtaining a set of $30$ instances in total. Each of these instances has at least one feasible solution ($x_i = \mu$ for all $i$). In Table~\ref{tab:eiqp} each line presents the results for one instance, $n$ indicates the size of the problem and $r$ the reference of the instance within the 5 instances with the same characteristics. The other columns are as described in Legends of Tables~\ref{tab:kcluster80} and~\ref{tab:bc140}--\ref{tab:ieqp}. We use the same experimental environment as described in Section~\ref{sec:kcluster}.

% \noindent \textbf{Experimental environment:} \\

% Experiments were carried out on a laptop with an Intel quad-core $i7$ processor of 1.73 GHz and $6$ GB of RAM using a Linux operating system. Again, we used the solver CSDP~\cite{csdp} for solving the semidefinite programs and the \textit{Conic Bundle} C-interface~\cite{cb} for the subgradient algorithm. We used the C-interface of \texttt{Cplex~12.5}~\cite{cplex125} for solving the quadratic programs. \\

\begin{table}
\centering
{\fontsize{6.5}{6.5}\selectfont
\begin{tabular}{c|c|c||c|c|c|c|c||c|c|c|c|c||c|c|c||c|c}
\multicolumn{3}{c||}{}&  \multicolumn{5}{|c||}{\texttt{MIQCR}} & \multicolumn{5}{|c||}{\texttt{MIQCR-CB $(p = 4{n \choose 2})$ }} & \multicolumn{3}{|c||}{\texttt{Cplex12.6}} & \multicolumn{2}{|c}{\texttt{Couenne}}   \\  
  \textit{class}& \textit{n}&\textit{r} & \textit{Gap}&  \textit{P1} &  \textit{P2} &  \textit{Tt}  & \textit{Nodes}& \textit{Gap}&  \textit{P1} &  \textit{P2} &  \textit{Tt}  & \textit{Nodes} &\textit{Gap}& \textit{Tt}  & \textit{Nodes}& \textit{Tt}  & \textit{Nodes} \\\hline
EIQP\_1&20&1&0.1&26&115&141&3956&0.1&2&37&39&1252 &159&2&2000& 16&2044\\
EIQP\_1&20&2&0.1&15&15&30&5&0.1&2&12&14&11&178&1&817&5&342\\
EIQP\_1&20&3&0.1&20&8&28&0&0.1&3&13&16&75&177&2&1000&16&1990\\
EIQP\_1&20&4&0.0&23&4&27&0&0.0&3&3&6&0&144&1&69&4&56\\
EIQP\_1&20&5&0.2&26&16&42&43&0.2&3&25&28&255&209&2&2000&8&475\\
EIQP\_1&30&1&0.0&289&108&397&80&0.0&14&45&59&142&182&25&8000&1644&116242\\
EIQP\_1&30&2&0.0&184&11&195&0&0.0&14&25&39&1&175&13&3500&347&16938\\
EIQP\_1&30&3&0.0&216&384&600&2741&0.0&37&59&96&218&183&33&13500&179&9714\\
EIQP\_1&30&4&0.1&218&68&286&23&0.1&19&23&42&26&160&218&53500&1231&53330\\
EIQP\_1&30&5&0.1&200&248&448&1251&0.1&20&92&112&1352&159&223&84000&596&33841\\
EIQP\_1&40&1&0.0&1093&2234&3327&1816&0.0&17&237&254&296&175&13&3000&148&1925\\
EIQP\_1&40&2&0.0&1219&452&1671&558&0.0&220&254&474&233&161&432&56000 & \multicolumn{2}{|c}{\texttt{failed} } \\
EIQP\_1&40&3&0.0&1325&196&1521&0&0.0&67&78&145&2&172&(10$\%$)&259500 & \multicolumn{2}{|c}{\texttt{failed}}\\
EIQP\_1&40&4&0.0&1672&870&2542&3226&0.0&1239&157&1396&799&169&2382&258000&(2$\%$)&87323\\
EIQP\_1&40&5&0.2&2281&863&3144&3871&0.2&1362&622&1984&6365 &163&(9$\%$)&294000&(14$\%$)&84268\\		
%\multicolumn{3}{c||}{\multirow{2}{*}{\textbf{Mean} }} & \multirow{2}{*}{\textbf{0.1}} & \multirow{2}{*}{\textbf{587}} & \multirow{2}{*}{\textbf{373}}& \multirow{2}{*}{\textbf{960}}& \multirow{2}{*}{\textbf{1171}} & \multirow{2}{*}{\textbf{0.1}}&\multirow{2}{*}{\textbf{201}}&\multirow{2}{*}{\textbf{112}}& \multirow{2}{*}{\textbf{314}}&  \multirow{2}{*}{\textbf{735}}&\multirow{2}{*}{\textbf{171}}&\textbf{349} & \multirow{2}{*}{\textbf{69259}}& \textbf{319}& \multirow{2}{*}{\textbf{34944}} \\
%\multicolumn{3}{c||}{ } & &&&&&&&&& &&\textbf{(13)} & & \textbf{(11)} & \\

\hline

EIQP\_2&20&1&0.1&23&66&89&2583&0.1&3&112&115&7979&159&2&2000&39&5585\\
EIQP\_2&20&2&0.0&60&10&70&3&0.0&10&12&22&11&180&3&3500&21&1748\\
EIQP\_2&20&3&0.0&30&17&47&73&0.0&2&18&20&140&138&1&500&5&336\\
EIQP\_2&20&4&0.2&32&54&86&1726&0.2&7&56&63&2777&155&2&3000&16&2132\\
EIQP\_2&20&5&0.1&19&31&50&341&0.1&6&45&51&1590&588&4&5000&9&724\\
EIQP\_2&30&1&0.3&481&494&975&4237&0.3&1296&165&1461&2011&186&86&50000&287&14978\\
EIQP\_2&30&2&0.1&337&256&593&1488&0.1&6&282&288&45378&171&21&13500&92&4655\\
EIQP\_2&30&3&0.0&264&201&465&1009&0.0&44&105&149&893&160&159&125000&284&15589\\
EIQP\_2&30&4&0.1&180&667&847&3728&0.1&6&226&232&4255&134&9&6000&86&4440\\
EIQP\_2&30&5&0.1&476&398&874&3591&0.1&1161&100&1261&1737&160&169&70000&413&24752\\
EIQP\_2&40&1&0.0&1078&1137&2215&739&0.0&33&750&783&2610&160&1739&236000&1387&39258\\
EIQP\_2&40&2&0.1&1216&2295&3511&4876&0.1&28&989&1017&5473&216&(11$\%$)&535500&(15$\%$)&82426\\
EIQP\_2&40&3&0.1&1118&2180&3298&3275&0.1& 17&1238&1255&6664&166&1366&253000&3016&95223\\
EIQP\_2&40&4&0.0&1306&165&1471&0&0.0&25&84&109&5&162&7&1500&80&391\\
EIQP\_2&40&5&0.0&2234&1415&3649&3616&0.0&67&697&764& 2841&177&(6$\%$)&450000& \multicolumn{2}{|c}{\texttt{failed}}	\\

%\multicolumn{3}{c||}{ \multirow{2}{*}{\textbf{Mean}} } &  \multirow{2}{*}{\textbf{0.1}}&  \multirow{2}{*}{\textbf{590}}&  \multirow{2}{*}{\textbf{626}} &  \multirow{2}{*}{\textbf{1216}} &  \multirow{2}{*}{\textbf{2086}}& \multirow{2}{*}{\textbf{0.1}}& \multirow{2}{*}{\textbf{181}}&  \multirow{2}{*}{\textbf{325}}&  \multirow{2}{*}{\textbf{506}}&  \multirow{2}{*}{\textbf{2902}}&  \multirow{2}{*}{\textbf{194}}& \textbf{373}&  \multirow{2}{*}{\textbf{116967}}& \textbf{582}&  \multirow{2}{*}{\textbf{23104}}\\
%multicolumn{3}{c||}{ } & &&&&&&&&&&& \textbf{(13)} & & \textbf{(13)} &

\end{tabular}
}
\caption {Computational results of \texttt{MIQCR-CB}, \texttt{MIQCR}, \texttt{Cplex12.6} and \texttt{Couenne} for 30 general integer instances with $20$, $30$ and $40$ variables of classes EIQP\textunderscore 1  and EIQP\textunderscore 2 (B\&B time limit: 3\;600 seconds).}
\label{tab:eiqp}

\end{table}
\par The numerical results comparing methods \texttt{MIQCR-CB}, \texttt{MIQCR}, \texttt{Cplex12.6} and \texttt{Couenne} are given in Table~\ref{tab:eiqp}. We observe that both algorithms \texttt{MIQCR-CB} and \texttt{MIQCR} solve all considered instances in less than $3\,649$ seconds of cpu time, while \texttt{Cplex12.6} solves only $26$ instances and \texttt{Couenne} only 24 instances, over the $30$ considered instances in less than one hour of cpu time. 
Comparing \texttt{MIQCR} and \texttt{MIQCR-CB}, we observe that \texttt{MIQCR-CB} is always faster than \texttt{MIQCR}. More precisely, while the initial gap remains the same for both methods, for class EIQP\textunderscore 1 (EIQP\textunderscore 2 resp.) the solution time for Phase~1, for \texttt{MIQCR-CB} in comparison to \texttt{MIQCR}, is divided on average by a factor of about $2.9$ ($3.3$ resp.), and by a factor $3.3$ ($1.9$ resp.) for the solution time of Phase~2. Hence the total run time significantly decreases (factor $3.1$ for EIQP\textunderscore 1 and $2.4$ for EIQP\textunderscore 2). 
Comparing \texttt{Cplex12.6} and \texttt{Couenne} with \texttt{MIQCR-CB}, we observe that \texttt{Cplex12.6} and \texttt{Couenne} are often faster than \texttt{MIQCR-CB} for the smallest considered instances, but with increasing dimension,  \texttt{MIQCR-CB} is superior to \texttt{Cplex12.6} and \texttt{Couenne}. Moreover, we can notice a significant decrease of the initial gap with method \texttt{MIQCR-CB} (by a factor $2600$) and of the number of visited nodes during the branch-and-bound algorithm (by a factor $180$ and $90$), in comparison to  \texttt{Cplex12.6} and \texttt{Couenne}. \\

\noindent \textbf{Experiments on the Integer Equipartition Problem} \\

We consider the Integer Equipartition Problem (IEP). This problem is an extension of the classical min-cut graph problem which consists in partitioning the vertices of a graph into a collection of disjoint sets satisfying specified size constraints, while minimizing the sum of weights of edges connecting vertices in different sets~\cite{HagKry99}. In (IEP) we consider $n$ types of items, $m$ items of each type, and a partition of the $n\cdot m$ items into $p$ equally sized sets. We assume that $n\cdot m$ is a multiple of $p$. For all pairs of type of items $(i,j), \,\, i \leq j$, we denote by $c_{ij}$ the cost of allocating each pair of items of types $i$ and $j$ to different sets. The problem consists thus to minimize the total cost of allocating the $n\cdot m$ items to the $p$ sets. By introducing decision variables $x_{ik}$ which represent the number of items of type $i$ allocated to set $k$,  (IEP) can be formulated as follows:

\begin{numcases}{(IEP)}
\min \, \, \displaystyle{\sum_{i<j}} \displaystyle{\sum_{k \neq l}} c_{ij} x_{ik} x_{jl} +\displaystyle{\sum_i} \displaystyle{\sum_{k<l}} c_{ii} x_{ik} x_{il}\nonumber \\ 
\mbox{s.t.} \nonumber \\
\quad \displaystyle{\sum_{i=1}^n}  x_{ik} = \frac{n\cdot m}{p} & $1\le k \le p$ \label{alloc-part}\\
\quad \displaystyle{\sum_{k=1}^p}  x_{ik} = m & $1\le i \le n$ \label{alloc-item} \\
\quad 0 \leq x_{ik} \leq \min (m,\frac{n\cdot m}{p}) & $1\le i \le n$,\,$1\le k \le p$ \\
\quad x_{ik} \in \mathbb{N}& $1\le i \le n$,\,$1\le k \le p$ 
\end{numcases}

\noindent where Constraints~(\ref{alloc-part}) ensure that exactly $\frac{n\cdot m}{p}$ items are allocated to each set, Constraints~(\ref{alloc-item}) ensure that all items are allocated to a set. This problem has a quadratic objective function, $n+p$ linear equalities and $n\cdot p$ general integer variables.

We generate instances of (IEP) with $n=3$ or $4$, $m=20$ or $24$ and $p=5$ or $6$, where the coefficients $c$ are uniformly distributed integers from the interval $[1,10]$. For each characteristics, we generate 10 instances. The numerical results comparing methods \texttt{MIQCR-CB}, \texttt{Cplex12.6} and \texttt{Couenne} are given in Table~\ref{tab:ieqp} where each line presents the results for one instance, $n$ indicates the number of types of items, $m$ the number of items of each type, $p$ the number of equally sized sets, and $r$ the reference of the instance within the 10 instances with the same characteristics. The other columns are as described in Legends of Tables~\ref{tab:kcluster80} and~\ref{tab:bc140}--\ref{tab:ieqp}

\begin{table}
\centering
\begin{scriptsize}
\begin{tabular}{c|c|c|c||c|c|c|c|c||c|c|c||c|c}
 \multicolumn{4}{c||}{}& \multicolumn{5}{|c||}{\texttt{MIQCR-CB $p=(4{n \choose 2})$}}& \multicolumn{3}{|c||}{\texttt{Cplex12.6}}& \multicolumn{2}{|c}{\texttt{Couenne}}   \\ 
  \textit{n}& \textit{m}& \textit{p}&\textit{r} & \textit{Gap}&  \textit{P1} &  \textit{P2} &  \textit{Tt}  & \textit{Nodes} &\textit{Gap}& \textit{Tt}  & \textit{Nodes}& \textit{Tt}  & \textit{Nodes}\\\hline
03&20&05&1&0.00&0&16&16&0&90.87&(5.07$\%$)&3054000&2888&1074696 \\
03&20&05&2&0.00&0&2&2&0&76.64&546&2244000&(1.4$\%$)&1271760  \\
03&20&05&3&2.07&1&110&111&29966&78.21&606&2487500&317&122622\\
03&20&05&4&1.50&1&38&39&8990&76.31&552&2244000&309&118438\\
03&20&05&5&1.73&0&13&13&2457&95.76&(2.17$\%$)&2875500&663&231750\\
03&20&05&6&0.11&1&13&14&1507&81.73&(8.24$\%$)&3272500&(2.1$\%$)&1182335\\
03&20&05&7&1.01&1&23&24&7481&97.8&(5.8$\%$)&3158000&2143&744678\\
03&20&05&8&0.11&9&8&17&571&81.6&(5.17$\%$)&3009000&1122&424302\\
03&20&05&9&0.12&1&13&14&1815&94.24&(4.73$\%$)&2951000&1534&595765\\
03&20&05&10&0.14&1&9&10&864&88.09&(5.07$\%$)&3054000&2194&814563\\
%\multicolumn{4}{c||}{\textbf{Mean} }&\textbf{0.7}&\textbf{2}&\textbf{25}&\textbf{27}&\textbf{5365}&\textbf{86.1}&\textbf{568 (3)}&\textbf{2834950}&\textbf{1396 (8)}& \textbf{658091}\\ 
\hline

03&24&06&1&0.00&1&5&6&0&88.2&(16.88$\%$)&2636000&(12.5$\%$)&661600\\
03&24&06&2&0.00&0&4&4&0&87.11&(18.6$\%$)&2575500&(14.1$\%$)&659827\\
03&24&06&3&0.00&1&7&8&29&80.19&(10.46$\%$)&2604500&(3.4$\%$)&721369\\
03&24&06&4&0.00&1&9&10&321&83.39&(11.17$\%$)&2501500&(3.8$\%$)&745568\\
03&24&06&5&0.00&5&3&8&0&84.66&(12.63$\%$)&2543500&(6.5$\%$)&713247\\
03&24&06&6&0.02&1&34&35&1943&88.75&(18.1$\%$)&2611500&(14.8$\%$)&650991\\
03&24&06&7&0.00&4&2&6&0&86.81&(15.61$\%$)&2503500&(10.8$\%$)&667285\\
03&24&06&8&0.35&20&241&261&31505&84.27&(14.99$\%$)&2582000&(11.5$\%$)&647872\\
03&24&06&9&0.04&11&7&18&309&77.35&(15.64$\%$)&2522000&(11.7$\%$)&693954\\
03&24&06&10&0.07&73&25&98&1754&79.58&(16.06$\%$)&2574000&(12.5$\%$)&662087\\
%\multicolumn{4}{c||}{\textbf{Mean} }&\textbf{0.1}&\textbf{12}&\textbf{34}&\textbf{46}&\textbf{3586}&\textbf{84}&\textbf{(15$\%$)}&\textbf{2565400}&\textbf{(10.1$\%$)}& \textbf{682380}\\ 
\hline

04&20&05&1&0.37&39&(0.18$\%$)&-&185288&94.02&(24.69$\%$)&1925000&(31.7$\%$)&546677\\
04&20&05&2&1.84&97&(1.26$\%$)&-&190884&91.46&(22.58$\%$)&1964500&(24.8$\%$)&571394\\
04&20&05&3&2.42&72&(1.83$\%$)&-&176167&86.75&(19.82$\%$)&2018500&(24.2$\%$)&590082\\
04&20&05&4&0.09&8&144&152&4714&90.17&(24.84$\%$)&1908500&(32.2$\%$)&528364\\
04&20&05&5&3.82&79&(2.33$\%$)&-&155941&90.05&(16.85$\%$)&2036500&(15.9$\%$)&630677\\
04&20&05&6&0.06&131&129&260&4454&84.43&(26.29$\%$)&1924500&(33.1$\%$)&518391\\
04&20&05&7&0.62&89&(0.46$\%$)&-&217944&87.1&(27.58$\%$)&1923000&(36.9$\%$)&496121\\
04&20&05&8&0.05&47&96&143&3942&97.29&(27.02$\%$)&1901500&(31.8$\%$)&512563\\
04&20&05&9&1.50&102&(1.15$\%$)&-&188395&95.93&(21.16$\%$)&199250&(24.6$\%$)&5764810\\
04&20&05&10&1.41&351&(1.27$\%$)&-&141855&91.96&(25.06$\%$)&1903000&(27.1$\%$)&557275\\
%\multicolumn{4}{c||}{\textbf{Mean} }&\textbf{1.2}&\textbf{102}&\textbf{123 (3)}&\textbf{185 (3)}&\textbf{126958}&\textbf{90.9}&\textbf{(23.6$\%$)}&\textbf{1949750}&\textbf{(28.2$\%$)}&\textbf{552803}\\ 
\hline

04&24&06&1&0.06&102&(0.06$\%$)&-&74597&93.62&(42.69$\%$)&1193500&(72.3$\%$)&308012\\
04&24&06&2&1.26&113&(1.12$\%$)&-&79635&90.36&(43.71$\%$)&1210500&(65.1$\%$)&314582\\
04&24&06&3&1.26&93&(1.09$\%$)&-&70921&93.48&(36.72$\%$)&1252500&(61.7$\%$)&336893\\
04&24&06&4&0.03&15&721&736&13603&99.26&(43.95$\%$)&1227000&(74.9$\%$)&298315\\
04&24&06&5&2.96&107&(2.29$\%$)&-&66534&94.02&(40.12$\%$)&1251000&(58.5$\%$)&340087\\
04&24&06&6&0.01&189&66&255&829&93.13&(45.61$\%$) &1224500&(74.2$\%$)&299641\\
04&24&06&7&0.00&136&353&489&6261&95.88&(43.85$\%$)&1199500&(71.5$\%$)&295719\\
04&24&06&8&0.25&16&(0.21$\%$)&-&118263&99.65&(48.48$\%$)&1218500&(75.0$\%$)&301299\\
04&24&06&9&0.01&151&260&411&3903&95.21&(42.23$\%$)&1194500&(60.3$\%$)&324965\\
04&24&06&10&0.85&1648&(0.81$\%$)&-&79072&92.82&(44.19$\%$)&1171500&(66.3$\%$)&315380\\
%\multicolumn{4}{c||}{\textbf{Mean} }&\textbf{0.7}&\textbf{257}&\textbf{350 (4)} &\textbf{378.4 (4)} &\textbf{51362}&\textbf{94.7}&\textbf{(43.2$\%$)}&\textbf{1214300}&\textbf{(68$\%$)}&\textbf{313489}
\\

\end{tabular}

\caption {Computational results of \texttt{MIQCR-CB}, \texttt{Cplex12.6}  and \texttt{Couenne} for 40 instances of (IEP) (time limit: 3\;600 seconds)}
\label{tab:ieqp}
\end{scriptsize}
\end{table}

We observe that \texttt{MIQCR-CB} is able to solve $27$ of the $40$ considered instances, while \texttt{Cplex12.6} solves only $3$ instances and \texttt{Couenne} 8 instances of the 10 of the smallest size. Moreover these 3 (8 resp.) instances are solved by \texttt{Cplex12.6} (\texttt{Couenne} resp.) within 568 (1396 resp.) seconds in average, while  \texttt{MIQCR-CB} solves the 10 instances within 27 seconds on average. An important advantage of \texttt{MIQCR-CB} is its gap that is on average 136 times smaller than the gap of \texttt{Cplex12.6}, and that is moreover quite stable independently from the characteristics of the instances.% As a consequence, the number of nodes visited during the branch-and-bound algorithm is much smaller for \texttt{MIQCR-CB} than for \texttt{Cplex12.6} and \texttt{Couenne} (factor 46 and 12 respectively).

\section{Conclusion}
\label{sec:conclusion}
We presented algorithm \texttt{MIQCR-CB} for solving general integer quadratic programs with linear constraints. This algorithm is an improvement of \texttt{MIQCR}~\cite{BEL12}. \texttt{MIQCR-CB} is an approach in two phases: the first phase calculates an equivalent quadratic reformulation of the initial problem by solving a semidefinite program, and the second phase solves the reformulated problem using standard mixed-integer quadratic programming solver. 

In Phase~1, a subgradient algorithm within a Lagrangian duality framework is used to solve the semidefinite program that yields the parameters for constructing an equivalent (convex) formulation. This significantly speeds up Phase~1 compared to the earlier method \texttt{MIQCR}. Furthermore, by construction, we can control the size of the reformulated problem. As a consequence, Phase~2 of our algorithm is also accelerated.  
Thus, \texttt{MIQCR-CB} builds a reformulated problem with a tight bound that can be computed in reasonable time, even for large instances.

 Computational experiments carried out on the $k$-cluster problem demonstrate that our method is competitive with the best current approaches devoted to the solution of this binary quadratic problems for instances with up to $120$ variables. Moreover, our approach is able to solve most of the instances with $160$ variables to optimality within 3~hours of cpu time. 
We also demonstrate that \texttt{MIQCR-CB} outperforms \texttt{MIQCR}, \texttt{Cplex12.6} and \texttt{Couenne} on two classes of general integer quadratic problems, which confirms the impact of the newly designed procedure for solving the SDP.

{\bf Acknowledgment:} We thank Franz Rendl for useful discussions and two anonymous referees for suggestions that improved this paper.

\bibliography{mybib}{}
\bibliographystyle{plain}

\end{document}